\documentclass[12pt]{amsart}
\usepackage{times}
\usepackage{amsfonts}
\usepackage{amsthm}
\usepackage{amsmath}
\usepackage{epic}
\usepackage{times}

\advance \topmargin by -\headheight
\advance \topmargin by -\headsep

\evensidemargin \oddsidemargin
\newtheorem{theorem}{Theorem}
\newtheorem{lemma}{Lemma}

\newtheorem{proposition}{Proposition}

\newtheorem{remark}{Remark}

\numberwithin{equation}{section}

\begin{document}  

\title{Rich, Sturmian, and Trapezoidal Words}

\author{Aldo de Luca}

\address{Dipartimento di Matematica e Applicazioni ``R. Caccioppoli", Universit\`a degli Studi di Napoli Federico II, Via Cintia, Monte S. Angelo, I-80126 Napoli, Italy}

\email{aldo.deluca@unina.it}

\author{Amy Glen}

\address{LaCIM, Universit\'e du Qu\'ebec \`a Montr\'eal, C.P. 8888, succursale Centre-ville, Montr\'eal, Qu\'ebec, H3C 3P8, Canada.}

\email{amy.glen@gmail.com}

\author{Luca Q. Zamboni}

\address{Institut Camille Jordan, 
Universit\'e Claude Bernard Lyon 1, 
43 boulevard du 11 novembre 1918, 
69622 Villeurbanne cedex
France}
\email{luca@unt.edu}


\date{February 2, 2008}

\begin{abstract}
In this paper we explore various interconnections between rich words, Sturmian words, and trapezoidal words.
Rich words, first introduced in  \cite{GJWZ} by second and third authors together with J. Justin and S. Widmer,  constitute a new class of finite and infinite words characterized by having the maximal number of palindromic factors. 
Every finite Sturmian word is rich, but not conversely \cite{DJP}.
Trapezoidal words were first introduced by the first author in studying the behavior of the subword complexity of finite Sturmian words.
Unfortunately this property does not characterize finite Sturmian words.
In this note we show that the only trapezoidal palindromes are Sturmian. More generally we show that Sturmian palindromes can be characterized either in terms of their subword complexity (the trapezoidal property) or in terms of their palindromic complexity. We also obtain a similar characterization of rich palindromes in terms of a relation between palindromic complexity and subword complexity.  \end{abstract}

\maketitle

\section{Introduction}

In \cite{DJP}, X. Droubay, J. Justin, and G. Pirillo showed that a finite word $W$ of length $|W|$ has at most $|W|+1$ many distinct palindromic factors, including the empty word. In
\cite{GJWZ}, the second and third authors together with J. Justin and S. Widmer initiated a unified study of both finite and infinite words characterized by this palindromic richness property. Accordingly we say that  a finite word $W$ is {\it rich} if and only if it has $|W|+1$ distinct palindromic factors, and an infinite word is rich if all of its factors are rich. Droubay, Justin and Pirillo showed that all Episturmian words (in particular all Sturmian words) are rich. Other examples of rich words are complementation symmetric sequences \cite{GJWZ}, symbolic codings of trajectories of symmetric interval exchange transformations \cite{FZ1, FZ2}, and certain $\beta$-expansions where $\beta$ is a simple Parry number \cite{AMPF}.

Let $u$ be a non-empty factor of a finite or infinite word $W.$ A factor of $W$ having exactly two occurrences of $u,$ one as a prefix, and one as a suffix, is called a {\it complete return} to $u$ in $W.$ In \cite{GJWZ}, the following fact is established:

\begin{proposition}\label{prop1} A finite or infinite word $W$ is rich if and only if for each non-empty palindromic factor $u$ of $W,$ every complete return to $u$ in $W$ is a palindrome.
\end{proposition}

In short, $W$ is rich if all complete returns to palindromes are palindromes. 
Given a finite or infinite word $W,$  let $C_W(n)$  (respectively $P_W(n))$ denote the {\it subword complexity function} (respectively the {\it palindromic complexity function}) which associates to each number $n\geq 0$ the number of distinct factors (respectively palindromic factors) of $W$ of length $n.$ Infinite Sturmian words are characterized by both their subword complexity and palindromic complexity. An infinite word $W$ is Sturmian if and only if $C_W(n)=n+1$ for each $n\geq 0.$ In \cite{DP}, X. Droubay and G. Pirillo showed that $W$ is Sturmian if and only if  $P_W(n)=1$ whenever $n$ is even, and $P_W(n)=2$ whenever $n$ is odd. In \cite{dL}, the first author studied the complexity function of finite words $W.$ He showed that if $W$ is a finite Sturmian word (meaning a factor of a Sturmian word), then the graph of
$C_W(n)$ as a function of $n$ (for $0\leq n\leq |W|)$ is that of a regular trapezoid: that is $C_W(n)$  increases by $1$ with each $n$ on some interval of length $r,$ then $C_W(n)$ is constant on some interval of length $s,$ and finally $C_W(n)$ decreases by $1$ with each $n$ on an interval of the same size $r.$ Such a word is said to be {\it trapezoidal}. More precisely,
for any word $W$ let us denote by  $R_W$ the smallest integer $p$ such that $W$ has no right  special
factor of length $p,$ and by $K_W$  the length of the shortest unrepeated
suffix  of $W$. Then we say that $W$ is a {\it trapezoidal word} if  and only if
 $|W| = R_W +K_W.$
However, in \cite{dL} the first author shows that the property of being trapezoidal does not characterize finite Sturmian words. For instance, the word $aaabab$ is not Sturmian although it is trapezoidal.\footnote{In \cite{dA}, F. D'Alessandro classified all non-Sturmian trapezoidal words.} 

The main results of this note are to give characterizations of both rich palindromes and Sturmian palindromes in terms of the palindromic complexity functions. We also show that every trapezoidal word is rich, but not conversely. In the case of rich palindromes we prove\footnote{An infinite version of Theorem~\ref{FGC} was obtained by the second and third authors together with M. Bucci and A. De Luca in \cite{BDGZ} using completely different methods.}:

\begin{theorem}\label{FGC} Let $W$ be a finite word. Then the following two conditions are equivalent:
\begin{itemize}
\item[(A)] $W$ is a rich palindrome.
\item[(B)] $P_W(n)+P_W(n+1)=C_W(n+1)-C_W(n) +2$ \,for each $0\leq n \leq |W|.$
\end{itemize} 
\end{theorem}

\vskip 2pt

\noindent While for Sturmian palindromes we prove\footnote{A different characterization of Sturmian palindromes was obtained by A. de Luca and A. De Luca in \cite{dLDL0}. See also \cite{dLDL1}.}:

\begin{theorem}\label{main} Let $W$ be a word of length $N.$ Then the following three conditions are equivalent:
\begin{enumerate}
\item[(A')] $W$ is a Sturmian palindrome.
\item[(B')] $P_W(n)+P_W(N-n)=2$ for each $0\leq n\leq N.$
\item[(C')] $W$ is a trapezoidal palindrome.
\end{enumerate}
\end{theorem}

\vskip 2pt

\section{Rich vs Trapezoidal Words}
In this section we show that all trapezoidal words are rich:

\begin{proposition}\label{prop2} Let $W$ be a trapezoidal word. Then $W$ is rich.
\end{proposition}

\begin{proof} We proceed by induction on $|W|.$  The result is clearly true if  $|W|\leq 2.$  Suppose every trapezoidal word of length less than $N$ is rich, and suppose that $W$ is trapezoidal (say on the letters $\{a,b\})$ of length $N.$ Let us suppose to the contrary that $W$ is not rich.
Then, by Proposition~\ref{prop1}, in $W$ there exists a complete return to some palindrome $P$ which is not a palindrome. Since, on a binary alphabet,  a complete return to a letter is always a palindrome, we can write (without loss of generality) that $P=aUa$ with $U$ possibly empty.
Since the prefix and suffix of  $W$ of length $N-1$ are both rich (by the induction hypothesis), it follows that $aUa$ is both a prefix and a suffix of $W,$ and that these are the only two occurrences of $aUa$ in $W.$  So $W$ itself is the complete return to $aUa$ which is not a palindrome. In particular $W$ is not a palindrome, which implies that $|W|\geq 2|aUa| +2.$

It follows that $K_W=|aUa|+1$  since $aUa$ occurs twice in $W$
and if some longer suffix of  $W$ occurred more than once in $W,$ then $aUa$ would occur at least three times in $W$.
Since $W$ is trapezoidal, we have  $R_W+K_W=|W|.$  Now the word W has a period $q = |W|-|aUa|= R_W+K_W - (K_W-1)= R_W+1.$ Let $\pi_W$ denote the minimal period of $W.$ 
Then $\pi_W \leq R_W+1.$
Since for any word $W,$  $\pi_W \geq R_W+1,$ it follows that $\pi_W = R_W+1.$
From Proposition 28 of \cite{dLDL0} we deduce that $W$ is a Sturmian,  and hence rich,  a contradiction.
\end{proof}

\begin{remark} We note that the converse is false; in fact $aabbaa$ is rich but not trapezoidal.
\end{remark}

\section{Proof of Theorem~\ref{FGC}}

We first show that (B) implies (A).
We assume $W$ satisfies (B). Taking $n=|W|$ and using $P_W(|W|+1)=C_W(|W|+1)=0$
and $C_W(|W|)=1,$ we deduce that $P_W(|W|)=1,$ and hence $W$ is a palindrome. It remains to show that $W$ is rich. Let $S$ denote the total number of distinct palindromic factors of $W.$ We will show that $S=|W|+1.$  Since $W$ itself is a palindrome we have 
\[S-1=\sum _{n=0}^{|W|-1}P_W(n)\]
Similarly since the empty word is a palindrome we have 
\[S-1=\sum_{n=1}^{|W|}P_W(n)\]
Thus
\begin{eqnarray*} 2S-2&=&\sum _{n=0}^{|W|-1}P_W(n) + \sum_{n=1}^{|W|}P_W(n)\\
&=&\sum_{n=0}^{|W|-1}(P_W(n) + P_W(n+1))\\
&=&\sum_{n=0}^{|W|-1}(C_W(n+1)-C_W(n) +2)\\
&=&C_W(|W|)-C_W(0) + 2|W|\\
&=&1-1 +2|W|=2|W|.
\end{eqnarray*}

\noindent Hence $S=|W|+1$ as required.

Next we show that (A) implies (B). We proceed by induction on the length of $W.$ The result is easily verified in the case $|W|\leq 2.$ Now suppose the result is true for all rich palindromes of length less than $N$ and suppose $W$ is a palindrome of length $N.$ Let $V$ denote the palindrome of length $N-2$ obtained by removing the first and last letter of $W.$ Since $V$ is also rich (see \cite{GJWZ}),  by the induction hypothesis we have 
$P_V(n)+P_V(n+1)=C_V(n+1)-C_V(n) +2$ \,for each $0\leq n \leq N-2.$

Let $N_0$ denote the length of a shortest factor $U$ of $W$ which is not a factor of $V.$ Then for $0\leq n< N_0-1$ we have  
$P_W(n)+P_W(n+1)=C_W(n+1)-C_W(n) +2.$

The word $U$ is either a prefix or a suffix of $W.$ We claim that it is in fact both a prefix and a suffix of $W,$ in other words a palindrome.
Suppose to the contrary that $U$ is not a palindrome. Without loss of generality we may assume that $U$ is a suffix of $W.$  Let $U'$ denote the longest palindromic suffix of $U.$ Since $|U'|<N_0,$ we have  $U'$ is also a factor of $V.$ Hence there exists a complete return $Z$ of $U'$ which is a proper suffix of $W.$ Since $W$ is rich, $Z$ is a palindrome. Since we are assuming that $U$ is not a palindrome and that $U'$ is the longest palindromic suffix of $U,$ it follows that $|Z|>|U|.$ Since $W$ is a palindrome, $Z$ is also a prefix of $W,$ and hence the proper suffix $U$ of $Z$ occurs in $V,$ a contradiction. Thus $U$ is a palindrome, and hence both a prefix and a suffix of $W.$ Thus $U$ is the only factor of $W$ of length $N_0$ which is not a factor of $V.$ Thus we have
\[P_W(N_0)=P_V(N_0)+1\,\,\,\,\mbox{and}\,\,\,C_W(N_0)=C_V(N_0)+1.\]

\noindent Since $P_V(N_0-1)+P_V(N_0)=C_V(N_0)-C_V(N_0-1)+2,$ $P_V(N_0-1)=P_W(N_0-1),$ and $C_V(N_0-1)=C_W(N_0-1),$  we deduce that
\[P_W(N_0-1) + (P_W(N_0)-1)=(C_W(N_0)-1)-C_W(N_0-1)+2\] and hence
\[P_W(N_0-1)+P_W(N_0)=C_W(N_0)-C_W(N_0-1)+2\]
in other words equality in (B) also holds for $n=N_0-1.$

We now claim that the only palindromic suffix of $W$ of length greater than $N_0$ is $W$ itself. In fact, if $W$ admitted a proper palindromic suffix of length greater than $N_0,$ then $U$ would be a factor of $V,$ a contradiction.
Thus we have 
\begin{eqnarray}P_W(n)=P_V(n)\,\,\,\mbox{for all}\,\, N_0<n<N.\label{P}\end{eqnarray}

Also, for each $N_0<n<N,$ let $UX$ (respectively $\bar{X}U)$ denote the prefix (respectively suffix) of $W$ of length $n,$ where $\bar{X}$ denotes the reversal of $X.$  Since $UX$ is not a palindrome it follows that $UX\neq \bar{X}U.$ Thus 
\begin{eqnarray}C_W(n)=C_V(n)+2 \,\,\,\mbox{for all} \,\, N_0<n<N.\label{C}\end{eqnarray}

We now verify (B) for $n=N_0.$ Starting with $P_V(N_0)+P_V(N_0+1)=C_V(N_0+1)-C_V(N_0)+2$ we obtain
\[(P_W(N_0)-1) + P_W(N_0+1)=(C_W(N_0+1)-2) -(C_W(N_0)-1) +2\]
\noindent and hence
\[P_W(N_0)+P_W(N_0+1)=C_W(N_0+1)-C_W(N_0)+2.\]

We next verify (B) for $N_0<n\leq N-2.$ Starting with $P_V(n)+P_V(n+1)=C_V(n+1)-C_V(n) +2$ and using (\ref{P}) and (\ref{C}) we obtain
\[P_W(n)+P_W(n+1)=(C_W(n+1)-2)-(C_W(n)-2)+2\]
\noindent and hence
\[P_W(n)+P_W(n+1)=C_W(n+1)-C_W(n)+2.\]

It remains to verify (B) for $n=N-1$ and $n=N.$ If $W$ is the constant word,
then  $P_W(N-1)=1,$ $P_W(N)=1,$  $P_W(N+1)=0,$ $C_W(N-1)=1,$ $C_W(N)=1,$ and $C_W(N+1)=0.$  Otherwise, 
$P_W(N-1)=0,$  $P_W(N)=1,$  $P_W(N+1)=0,$ $C_W(N-1)=2,$ $C_W(N)=1,$ and $C_W(N+1)=0.$ In either case one readily verifies
(B) for $n=N-1$ and $n=N.$ 
This completes the proof of Theorem~\ref{FGC}.

\section{Proof of Theorem~\ref{main}}

We begin with the following lemma:

\begin{lemma} Let $W$ be a word of length $N$ satisfying either condition of Theorem~\ref{main}. Then $W$ is a rich palindrome. Hence by Theorem~\ref{FGC} we have $P_W(n)+P_W(n+1)=C_W(n+1)-C_W(n)+2$ for $0\leq n\leq N.$
\end{lemma}

\begin{proof} Since any Sturmian word is trapezoidal, by Proposition~\ref{prop2} one has that if $W$ satisfies either condition (A') or (C'), then it is rich. Let us suppose that $W$ satisfies condition (B').  Since $P_W(N)=P_W(0)=1,$ we have  $W$ is a palindrome.
To see that $W$ is rich, let $S=P_W(0)+P_W(1)+P_W(2)+...+P_W(N)$ denote the number of distinct palindromic factors of $W.$
Then 

\begin{eqnarray*} 2S&=& P_W(0)+P_W(N) + P_W(1)+P_W(N-1)+\ldots +P_W(N) + P_W(0)\\&=&2(N+1).\end{eqnarray*}
\noindent Whence $S=N+1=|W|+1.$

\end{proof}

 We note that condition (B') is equivalent to saying that the word $P_W(0)P_W(1)P_W(2)...P_W(N)$  is a $\theta$-palindrome on the alphabet $\{0,1,2\}$
with respect to the involutory antimorphim $\theta$ defined by $\theta (0)=2,$ $\theta (2)=0$ and $\theta (1)=1.$

Assume first that $W$ is a Sturmian palindrome. For $0\leq n\leq N-1,$ set $D_W(n)=C_W(n+1)-C_W(n).$  In \cite{dL}, the first author showed that the word $D_W(0)D_W(1)D_W(2)....D_W(N-1)$ is of the form $1^r0^s(-1)^r.$  In other words, that $W$ is a {\it trapezoidal word:} $C_W(n)$ increases by $1$ with each $n$ on an interval of length $r,$ then stabilizes, and eventually decreases by $1$ with each $n$ on an interval of the same size $r.$  
The trapezoidal property of $W$ together with the preceding lemma imply that the word $P_W(0)P_W(1)P_W(2)...P_W(N)$ begins with a block of the form $121212\ldots$ (corresponding to the interval of length $r$ on which $C_W(n+1)-C_W(n)=1)$, and terminates with a block of the form $\ldots 010101$ (corresponding to the interval on which $C_W(n+1)-C_W(n)=-1)$,  and moreover by the trapezoidal property, these two blocks are of the same length. Between these two blocks is either a block of the form $11\ldots 11$ or of the form $202\ldots 020$ corresponding to the interval on which 
$C_W(n+1)-C_W(n)=0.$ Hence $W$ satisfies condition (B').

Next suppose $W$ satisfies (B').  First observe that for each $n$ we have $P_W(n)\in \{0,1,2\},$ and $P_W(1)\neq 0.$  
If  $P_W(1)=1,$ then
$W$ is equal to the constant word, and hence a Sturmian palindrome.  Next suppose $P_W(1)=2.$ In this case $W$ is a binary palindromic word, say on the alphabet $\{a,b\}.$
To show that $W$ is Sturmian, it suffices to show that $W$ is {\it balanced,} i.e., given any two factors $u$ and $v$ of $W$ of the same length, we have
$||u|_a-|v|_a|\leq 1,$ where $|u|_a$ denotes the number of occurrences of the letter $a$ in $u.$ Suppose to the contrary that $W$ is not balanced. Then, it is well known (see for instance  Proposition 2.1.3 in \cite{Lo}) that there exists a palindrome $U$ such that both $aUa$ and $bUb$ are factors of $W.$ Thus $W$ contains two distinct palindromes of the same length, which implies that $|U|$ is odd. For otherwise, if $|U|$ were even, then taking
$k=2^{-1}|U| +1,$ we have $P_W(2k)=2,$ and hence by (B') $P_W(N-2k)=0,$ and
hence $P_W(N)=0,$ a contradiction. Since $W$ is a palindrome and contains both $aUa$ and $bUb,$ the palindrome $U$ must have at least two complete
returns in $W,$ one beginning in $Ua,$  which we denote by $X,$ and one beginning in $Ub,$ which we denote by $Y.$ Since $W$ is rich we have  both $X$ and $Y$ are palindromes with $X\neq Y.$ 

If both $|X|$ and $|Y|$ are greater than $|U|+1,$ then both $|X|$ and $|Y|$ must be even. In fact, suppose to the contrary that $|X|$ were odd. Then $|X|\geq |U|+2.$ But then $W$ would contain three palindromes of length $|U|+2,$ namely $aUa,$ $bUb,$ and the central palindromic factor of length $|U|+2$ of $X$ which is necessarily distinct from
both $aUa$ and $bUb$ since $X$ cannot contain an occurrence of $U$ other than as a prefix and as a suffix. The same argument shows that $|Y|$ must be even. Without loss of generality we can assume $|X|\leq |Y|.$ Then, as $X$ and the central palindrome of $Y$ of length $|X|$ are distinct, it follows that
$W$ contains two distinct palindromes of even length $|X|.$ Thus, $P_W(|X|)=2,$ and hence $P_W(N-|X|)=0,$ and hence $P_W(N)=0,$ a contradiction.

Thus it remains to consider the case in which either $|X|$ or $|Y|$ is equal to $|U|+1.$ Without loss of generality suppose $|X|=|U|+1.$ This means that
$X=Ua=aU$ and hence $U$ is the constant word $U=a^{|U|}.$
In this case $|Y|\geq |U|+2$ and by the previous argument must be even.
But then $X$ and the central palindrome of $Y$ of length $|X|$ are two distinct palindromic factors of $W$ of even length, a contradiction.
Thus we have shown that conditions (A') and (B') are equivalent.

Now we show that (A') is equivalent to (C'). The first author showed in \cite{dL} that every finite Sturmian word is trapezoidal. Thus (A') implies (C'). 
To see that (C') implies (A'), we proceed by induction on $|W|.$  The result is clearly true if $|W|\leq 2.$  Next suppose the result is true for $|W|<N $ and let $W$ be a trapezoidal palindrome of length $N.$ Since a trapezoidal word is necessarily on a two-letter alphabet, say $\{a,b\},$ we can write, without loss of generality, $W=aVa.$ Then $V$ is a trapezoidal palindrome, since factors of trapezoidal words are trapezoidal (see \cite{dA}).  By the induction hypothesis, $V$ is a Sturmian palindrome. 
If $W$ is not Sturmian, then there exists a palindrome $U$ such that $aUa$ and $bUb$ are factors of $W.$
Since $V$ is Sturmian, we have  $aUa$ is both a prefix and suffix of $W,$ and $bUb $ is a factor of $V.$  Since in $V,$ all complete returns to $U$ are palindromes, between an occurrence of $bUb$ in $V$ and the suffix $aU$ of $V$ there must be an occurrence of $bUa.$
Since $V$ is  a palindrome we have  $aUb$ is also a factor of $V.$
Hence each of $aUa, bUb, aUb,$ and $bUa$ is a factor of $W.$ This implies that both $aU$ and $bU$ are right special factors of $W,$ a contradiction since the trapezoidal property implies that for any $0\leq n\leq |W|,$ there exists at most one right special factor of $W$ of length $n.$ Thus $W$ must be Sturmian.
This concludes our proof of Theorem~\ref{main}.
 
\noindent {\bf Remark:}  A. De Luca \cite{DL} suggested the following alternate simple proof that (C') implies (A'): Let $W$ be a trapezoidal palindrome. Without loss of generality we can assume that $|W|\geq 2,$ for otherwise the result is clear.
Let $U$ denote the longest proper palindromic suffix of $W.$ Since $W$ is a palindrome,  $U$ is the longest border of $W,$ whence $|W|=\pi_W +|U|.$ By Proposition~\ref{prop2}, $W$ is rich, hence $U$ is the longest repeated suffix of $W.$ Thus $K_W=|U|+1.$ Since $W$ is trapezoidal we have that
$\pi_W=|W|-|U|=R_W+K_W-|U|=R_W+1.$ By Proposition 28
of  \cite{dLDL0} we deduce that $W$ is Sturmian.

\end{document}